%% file: FinalManuscript.tex
\documentclass[twocolumn]{autart}

\include{pkg_hdr}
\newif\iftext
\texttrue

\newif\ifproof
\prooftrue

\newif\ifhighlightthm
\highlightthmfalse

\ifhighlightthm
	\definecolor{theorem_highligts}{RGB}{255,127,0}
\else
	\definecolor{theorem_highligts}{RGB}{0,0,0}
\fi

\ifproof
   
\else
  \excludecomment{inproofs}
\fi

\excludecomment{toremove}
\newtheorem{theorem}{{\color{theorem_highligts} Theorem}}[section]
\newtheorem{definition}{{\color{theorem_highligts} Definition}}
\newtheorem{lemma}[theorem]{{\color{theorem_highligts} Lemma}}
\newtheorem{proposition}[theorem]{{\color{theorem_highligts} Proposition}}
\newtheorem{corollary}[theorem]{ {\color{theorem_highligts} Corollary}}

\newcommand{\mset}{\mathbf{\Sigma}}
\newcommand{\graph}{\mathbf{G}}
\newcommand{\reels}{\mathbb{R}}
\newcommand{\tran}{\top}
\newcommand{\mulno}{\mathcal{H}} 

\newcommand{\kron}{\otimes}
\newcommand{\jsr}{\hat \rho}
\newcommand{\card}[1]{N_{#1}}
\newcommand{\bigO}{\mathcal{O}}

\begin{document}

\begin{frontmatter}

\title{Stability of discrete-time switching systems with constrained switching sequences.} 

\author[LLN]{Matthew Philippe\thanksref{fnrs}}\ead{matthew.philippe@uclouvain.be},    %
\author[Illi]{Ray Essick\thanksref{illi}}\ead{ressick2@illinois.edu},
\author[Illi]{Geir Dullerud\thanksref{illi}}
\ead{dullerud@illinois.edu},
\author[LLN]{Rapha\"{e}l M. Jungers\thanksref{fnrs}}\ead{raphael.jungers@uclouvain.be}.

\address[LLN]{Institute of Information and Communication Technologies, Electronics and Applied Mathematics, Department of Mathematical Ingeneering (ICTEAM/INMA) at the Universit\'e
 catholique de Louvain, B-1348 Louvain-la-Neuve, Belgium}  
\address[Illi]{Department of Mechanical Science and Engineering and the Coordinated Science Laboratory, University of Illinois at Urbana-Champaign. Urbana, IL 61801, USA.}  
\thanks[fnrs]{M. Philippe is a F.N.R.S./F.R.I.A. fellow; R. Jungers  is a F.R.S./F.N.R.S. research associate.
They are supported by the Belgian Interuniversity Attraction Poles, and by the ARC grant 13/18-054 (Communaut\'{e} fran\c{c}aise de Belgique).
}
\thanks[illi]
{R. Essick and G.E. Dullerud were partially supported by
grants NSA SoS W911NSF-13-0086 and AFOSR MURI FA9550-10-
1-0573}
          

\begin{abstract}                          
We introduce a novel framework for the stability analysis of discrete-time linear switching systems with switching sequences constrained by an automaton. 
The key element of the framework is the algebraic concept of multinorm, which associates a different norm per node of the automaton, 
and allows to exactly characterize stability. 
Building upon this tool, we develop the first arbitrarily accurate approximation schemes for estimating the \emph{constrained} joint spectral radius $\jsr$, 
that is the exponential growth rate  of a switching system with constrained switching sequences. 
More precisely, given a relative accuracy $r > 0$, the algorithms compute  an estimate of $\jsr$ within the range $[\jsr, (1+r)\jsr]$.
 These algorithms amount to solve a well defined convex optimization program with known time-complexity, and whose size depends on the desired relative accuracy $r>0$.  
\end{abstract}

\end{frontmatter}

\section{Introduction}

\begin{figure}[ht]
\centering
\begin{subfigure}[b]{0.2\textwidth}
\centering
\includegraphics[scale = 0.4]{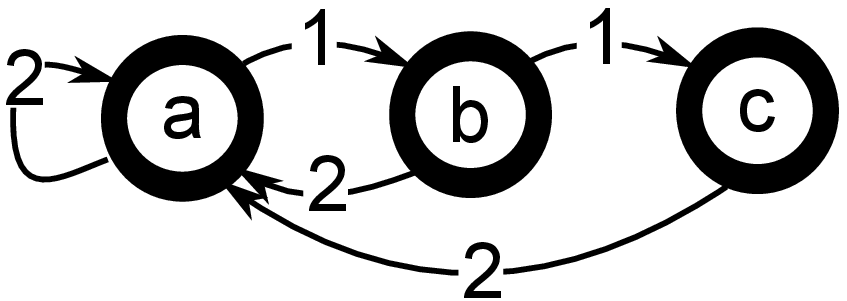}
\caption{}
\label{fig:intro_exa}
\end{subfigure}
~
\begin{subfigure}[b]{0.2\textwidth}
\centering
\includegraphics[scale = 0.4]{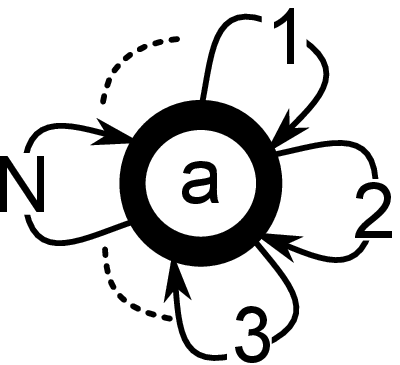}
\caption{}
\label{fig:intro_exb}
\end{subfigure}
\caption{The labels are represented on the edges.  Fig. \ref{fig:intro_exa} corresponds to the example in the introduction, where mode ``1'' cannot occur more than twice in a row. Node ``a'' is reached when the controller works, ``b'' is reached after one failure, and ``c'' after two failures.  The automaton of Fig. \ref{fig:intro_exb} accepts arbitrary switching sequences on $N$ modes.}
\label{fig:intro_ex}
\end{figure}

In this paper, we study  discrete-time linear switching systems having the particularity that their switching sequences are constrained by logical rules. We begin with an example introducing such systems.\\*
Given an \emph{unstable} matrix $A_1 \in \reels^{n \times n}$ and an input-to-state matrix
$B \in \reels^{n \times m}$, one computes a control gain matrix $K \in \reels^{m \times n}$
such that $A_2 = (A_1+BK)$ is stable.
The matrix $A_2$ dictates the closed-loop dynamics of a plant, $x_{t+1} = A_2 x_t$, whose stability is ensured by a state-feedback controller.
Let us now consider that the controller can fail at any time $t$, such that the dynamics at that time are given by
$ x_{t+1} = A_1 x_t. $
Then, the dynamics of the plant with failures can be modelled as a switching system $$
 x_{t+1} = A_{\sigma(t)} x_t, \vspace{-0.25cm}$$
where $\sigma(t) \in \{1,2\}$ is the \emph{mode} of the system and 
$\sigma(0),\,  \sigma(1), \, \ldots,$ is the \emph{switching sequence} that drives the system.
Without more information on the occurrences of the failures, we can only assume that the system is unstable.
Indeed, in the case of a permanent failure, represented by the switching sequence $ \sigma(t) = 1, \, \forall t \geq 0,$ the plant would follow the unstable dynamics 
$x_{t+1} = A_1 x_t$ at every time $t \geq 0$.
However, if we knew with \emph{certainty} that the failure cannot occur \emph{more than twice} in a row, 
then the above switching sequence would no longer be possible, and the system could very well be stable. \\*
This paper provides tools for the stability analysis of switching systems with constrained switching sequences, as in the example above. We say that the switching system on the matrix set $\mset = \{A_1, A_2, \ldots, A_N\}$ is stable if and only if,  
for all \emph{accepted} switching sequences $\sigma(0), \sigma(1), \ldots$, we have $ \lim_{t \rightarrow \infty} A_{\sigma(t)} \cdots A_{\sigma(0)} = 0.$   \\*
 Switching systems find applications in many theoretical and engineering related domains \cite{JuHeCOLS, HeMiOAMS, JuDIFSOD, JuTJSR, LiMoBPIS, OlReCPIN}, and the stability of switching systems  is known to be a challenging question \cite{ShWiSCFS, LiAnSASO, LiMoBPIS}. \\*
If one does not impose any constraint on switching sequences, the resulting system is called an \emph{arbitrary} switching system.
These systems have received a lot of attention in the past (e.g.,
\cite{AgLiLASC, AhJuJSRA, AnShSC, JuTJSR}).
The stability of an arbitrary switching system on a set of matrices $\mset$ is characterized by its \emph{joint spectral radius} (JSR) $\jsr(\mset)$ (introduced in \cite{RoStANOT}). It represents the worst case exponential growth rate of the system, and stability is equivalent to  $\jsr(\mset) < 1$, which is also equivalent to exponential stability. 
There has been a lot of research effort towards the computation and approximation of the JSR (see e.g. \cite{AhJuJSRA, JuTJSR, BlNeOTAO, VaHeJSRA} and references therein). 
One common way to do so is by computing a contractive invariant norm for the system \cite{AnShSC, BlNeOTAO, PaJaAOTJ, AnLaASCF }, which \emph{always exists} for stable arbitrary switching systems. For any level of relative accuracy $r > 0$, one can approximate these norms with quadratic/sum-of-square polynomials  \cite{BlNeOTAO, PaJaAOTJ} and
provide an \textit{upper bound}
 on the joint spectral radius within the range $[\jsr, (1+r) \jsr]$. 
 The computation of this estimate is done with finite time-complexity.  \\*
Our focus is on the stability of switching systems having logical rules on their switching sequences, 
such as the ones studied in \cite{AhJuJSRA, BlFeSAOD, DaAGTS, EsLeCOLS, KoTBWF, KuChSSSF, LeDuODAF, LeDuUSOD, LiAnSASO, PhJuASCF, WaRoSOLA}. 
We refer to these as \emph{constrained switching systems}, and represent the rules by using an \emph{automaton}. 
An automaton is a strongly connected, directed and 
labelled graph $\graph(V,E)$,
with $\card{V}$ nodes in $V$ and $\card{E}$ edges in $E$. 
The edge $(v, w, \sigma) \in E$ between the two nodes $v,w \in V$ carries the \emph{label} $\sigma \in \{1, \ldots, N\}$, which maps to a mode of the switching system. 
A sequence of modes $\sigma(0), \sigma(1), \ldots,$  is \emph{accepted} by the graph $\graph$ if there is a path in $\graph$ carrying the sequence as the succession of the labels on its edges. 
We do not specify an initial and final node for accepted paths, 
in that we depart from the usual definition for an automaton (see \cite{LoACOW}, Section 1.3). The accepted switching sequences form a symbolic dynamical system called \emph{sofic shift} (see \cite{LoACOW}, Section 1.5).
Examples of automata are given in Figure \ref{fig:intro_ex}. 

The system on the automaton $\graph$ with matrix set $\mset$ is denoted $S(\graph, \mset)$. The stability of $S(\graph, \mset)$ is characterized 
by the \emph{ constrained joint spectral radius}, introduced by Dai \cite{DaAGTS}.  A proof of the following is given in Annex \ref{annex:proofCJSR}.
\begin{theorem}[Dai \cite{DaAGTS}, Corollary 2.8] 
A  \mbox{constrained} switching system $S(\graph,\mset)$ 
is stable if and only if its \emph{constrained joint spectral radius} (CJSR),
defined as 
\begin{equation}
\begin{aligned}
& \jsr(S) \overset{\Delta}{=} \lim_{t \rightarrow \infty} \max_{\sigma({\cdot})} \{ \| A_{\sigma(t-1)} \cdots A_{\sigma(0)} \|^{1/t}: \\
& \qquad \sigma(0), \ldots, \sigma(t-1) \text{ is accepted by  $\graph$} \},
\end{aligned}
\label{eq:jsr}
\end{equation}
 satisfies $\jsr(S) < 1$. This also implies exponential stability.
For all accepted switching sequences,
$$\exists K \geq 1, 0 < \rho < 1: \forall T \geq 0, \|A_{\sigma(T-1)} \cdots A_{\sigma(0)}\| \leq K \rho^T.$$ 
\label{thm:stability-cjsr}
\end{theorem}

The CJSR, defined as (\ref{eq:jsr}), is independent of the norm used and homogeneous in $\mset$.\\
To the best of our knowledge, previous works on the stability of constrained switching systems have focused on establishing algorithmically checkable stability conditions, without studying their conservativeness.
  There is a particular interest  in using multiple quadratic Lyapunov functions as stability certificates \cite{BlFeSAOD, LeDuUSOD,  LeKhDASO,  EsLeCOLS, BrMLFA, LiAnSASO}. These approaches provide sets of LMIs whose feasibility is sufficient for stability. In \cite{BlFeSAOD, LeDuUSOD}, a hierarchy of more and more complex LMIs is presented such that, for any stable system, all LMIs starting from a certain level of complexity (depending on the system) are feasible.
 The methods discussed above can be used to obtain \emph{upper bounds} 
 on the CJSR (a feasible LMI indicates  $\jsr(S) < 1$). However, accuracy guarantees on these bounds, similar to that existing on the JSR estimation, have not been proven yet. 
\\*  
The framework we introduce allows to obtain accuracy guarantees\footnote{Preliminary results were presented in \cite{PhJuCLTF}.}.   A direct approach could rely on building an arbitrary switching system whose JSR equals the CJSR of the constrained system \cite{KoTBWF, WaRoSOLA}.
  We provide more efficient and intuitive techniques.  We generalize the recent results from \cite{AhJuJSRA} towards constrained switching systems. In \cite{AhJuJSRA}, the authors focus on systems $S(\graph, \mset)$	with $\graph$ accepting arbitrary switching sequences, and provide accuracy bounds for the JSR estimation using multiple Lyapunov functions. The generalization of these results to general constrained switching systems was left as an open question.
\\*
The plan of the paper is as follows.
Section \ref{sec:mn} introduces the algebraic concept of multinorm, which characterizes the stability of constrained switching systems
as contractive norms do for arbitrary switching systems.
  In Section \ref{sec:approx}, we focus on the more algorithmic question of the approximation - in finite time and with arbitrary accuracy - of the CJSR of a system $S(\graph, \mset)$.
In Section \ref{sec:example}, we illustrate our framework on a numerical example.\\*
\textbf{Notations.} 
The matrix $A^\tran \in \reels^{m \times n}$ is the transpose of  $A \in \reels^{n \times m}$. 
A path $p$ of length $T \geq 0$  in a graph $\graph$ is a sequence of $T$ consecutive edges. For a path $p$ with length $T \geq 1$, by a slight abuse of notations, we let 
$A_p = A_{\sigma(T)} \cdots A_{\sigma(1)}, $
where the $\sigma(1), \ldots, \sigma(T)$ are the $T$ labels along $p$. If $T = 0$, we let $A_p = I$, the identity matrix of $\reels^n$.

\section{Lyapunov functions for constrained switching systems}
\label{sec:mn}
The stability of \emph{arbitrary} switching systems is equivalent to the existence of a contractive norm serving as a Lyapunov function. We recall that a norm is a sub additive, positive definite  and  homogeneous function.
\begin{proposition}[e.g. \cite{JuTJSR}, Proposition 1.4]
The joint spectral radius of a set of matrices $\mset$ is given by 
\begin{equation}
 \jsr(\mset) \overset{\Delta}{=} \inf_ {|\cdot|}  \min_{ \gamma } 
 \Big \{ \gamma: \,µ |Ax| \leq \gamma |x|, \, \forall x \in \reels^n, \,  A \in \mset \Big \} .
 \label{eq:JSRnorm2}
\end{equation}
where the infimum is taken over all vector norms in $\reels^{n}.$
\label{prop:JSRnorms}
\end{proposition}
A stable arbitrary switching system has $\jsr(\mset) < 1$ (see \cite{JuTJSR}) and from Proposition  $\ref{prop:JSRnorms}$, there exists a norm $|{\cdot}|$ such that
$|Ax| < |x|$, $\forall x \in \reels^n$, $\forall A \in \mset$. 
It is however straightforward to build stable \textit{constrained} switching systems 
for which contractive norms \emph{do not} exist. 
\begin{exmp}
A scalar arbitrary switching  system built on $\mset = \{A_1, A_2\} = \{2, 1/8\}$ is unstable and it has no contractive norm due to $A_1$.
Consider now the  automaton $\graph$ of Figure \ref{fig:scalar12aut}.
\begin{figure}[ht]
\centering
\includegraphics[scale = 0.4]{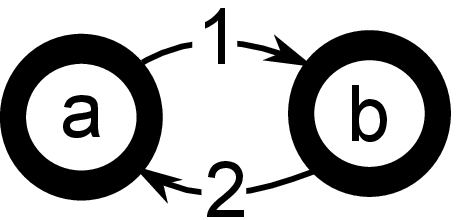}
\caption{}
\label{fig:scalar12aut}
\end{figure}
The periodic system $S(\graph, \mset)$ is  stable and, applying (\ref{eq:jsr}),
$\jsr(S) = \lim_{t \rightarrow \infty}\|(A_1A_2)^t\|^{1/(2t)} = 1/2$. 

\label{exmp:scalarswitching}
\end{exmp}
We fill this gap between arbitrary switching and constrained switching systems by introducing the algebraic concept of \emph{multinorm}.
\begin{definition}[Multinorm]
A \emph{multinorm} $\mulno$ for a system $S(\graph(V,E), \mset)$ is a set of
$\card{V}$ norms $\mulno = \{|{\cdot}|_v, \, v \in V\}$.
The \emph{value} $\gamma^*(\mulno)$ of a multinorm is defined as
\begin{equation} 
\begin{aligned}
&\gamma^*(\mulno) \overset{\Delta}{=} \min_\gamma \{ \gamma : |A_\sigma x|_w \leq \gamma |x|_v,\\
& \qquad \qquad \forall x \in \reels^n, (v,w,\sigma) \in E \}.
\end{aligned}
\label{eq:mn}
\end{equation}

\label{def:mn1}
\end{definition}
 Similar ideas have appeared in the literature \cite{LeDuUSOD, AhJuJSRA, BrMLFA, AhARAH, DaRiSAAC},
  where  multiple Lyapunov functions are considered for characterizing stability where single Lyapunov functions fail to do so. 
 Their role was either to provide a sufficient stability condition under the form of a set of LMIs, or in
\cite{AhJuJSRA}, to characterize the stability of \emph{arbitrary} switching systems using multiple Lyapunov functions. In comparison, we provide general \emph{necessary and sufficient} conditions for the stability of constrained switching systems, using multiple Lyapunov functions with exactly one \emph{norm} per \emph{node} of $\graph$.
\begin{proposition}
The \emph{constrained joint spectral radius} (\ref{eq:jsr}) of a system $S(\graph, \mset)$ satisfies 
\begin{equation}
\hat \rho(S) \overset{\Delta}{=} \inf_{\mulno} \{ \gamma^*(\mulno): \mulno \text{ is a multinorm for $S$} \}.
\end{equation}
\label{prop:CJSRmn}
\end{proposition}\vspace{-1cm}

\begin{pf}
We first show that the value of any multinorm for a system is an upper bound of its CJSR.\\*
Consider a multinorm $\mulno = \{|{\cdot}|_v, \, v \in V\}$ for $S$ with value $\gamma$. For any path $p$ with length $k \geq 1$ between two nodes $v$ and $w$ in $\graph$, from (\ref{eq:mn}), we get 
$ |A_px|_w \leq \gamma^k |x|_v.$\\*
For any norm $|{\cdot}|$, by equivalence of norms in $\reels^n$, there exists $0 < \alpha < \beta$ such that the inequalities
 $  \alpha |x| \leq |x|_v \leq \beta |x| $
 hold for all $x \in  \reels^n$ and all the norms $|{\cdot}|_v$ in $\mulno$.\\*
Considering  the classical definition for an induced matrix norm, we have
\begin{equation*}
\begin{aligned}
 \|A_{p}\| & \overset{\Delta}{=}  \max_{|x| = 1}\frac{|A_{p}x|}{|x|}
  \leq \frac{\beta}{\alpha}\max_{|x| = 1}\frac{|A_{p}x|_w}{|x|_v}
  \leq \frac{\beta}{\alpha} \gamma^k.
 \end{aligned}
 \label{eq:boundedness}
\end{equation*}
Taking paths of lengths $k \rightarrow \infty$ and the $k$th root of the above inequality, we obtain  $\jsr(S) \leq \gamma$ from (\ref{eq:jsr}).\\*
We now show that for any $\epsilon > 0$ there exists a multinorm of value at most $(\jsr(S) + \epsilon)$.
Consider the scaled set of matrices  
$$ \mset' = \{A'_i = A_i /  ( \jsr(S) + \epsilon ), \, i = 1, \ldots, N \}.$$  
The CJSR of $S(\graph, \mset)$ is an homogeneous function of $\mset$ (see (\ref{eq:jsr})). The system $S(\graph, \mset')$ is then stable since $\jsr(S)/( \jsr(S) + \epsilon) < 1$.\\
We  define, at each node $v \in V$, the following functions which we then prove are norms:
\begin{equation}
\begin{aligned}
&|x|_v := \sup_p\{ |{A}_p' x| : 
 p \text{ is a path with origin $v$} \} ,
\label{eq:nodconst}
\end{aligned}
\end{equation} 
where $|{\cdot}|$ is the euclidean norm.
These functions are well-defined 
(by exponential stability of $S(\graph, \mset')$), 
sub-additive, homogeneous, 
and positive definite (with paths of length 0, $|x|_v \geq |x|$), hence they are norms.
Moreover, for any edge $(v, w, \sigma) \in E$, and all $x \in \reels^n$, we have
\begin{align}
	|x|_v &= \sup_p\{ |{A}_p' x| : p \text{ is a path with origin $v$} \}, \nonumber\\
	&\geq \sup_q\{ |{A}_q'A'_\sigma x| : q \text{ is a path with origin $w$} \}, \label{eq:middleineq}\\
	 &= |A'_\sigma x|_w, \nonumber
\end{align}
where (\ref{eq:middleineq}) is obtained by taking $p$ starting with the edge $(v,w,\sigma)$.
Since $ A_{\sigma}' = A_{\sigma} / ( \jsr(S) + \epsilon)$, we have
$| A_\sigma x |_{w} \leq (\jsr(S)+\epsilon) | x |_{v}$
for all $(v, w, \sigma) \in E$. Thus, the value of the multinorm with the norms defined in (\ref{eq:nodconst}) is upper bounded by $\jsr(S) + \epsilon$.
\end{pf}
We conclude that a constrained switching system is stable if and only if it has a multiple Lyapunov function with exactly one norm per node of its graph. The proof of the following is direct from Proposition \ref{prop:CJSRmn}.
 \begin{theorem}
A constrained switching system $S( \graph, \mset)$ is stable if and only if 
it admits a multinorm $\mulno$ with value $\gamma^*(\mulno) < 1$.
Such a multinorm is said to be a \emph{Lyapunov} multinorm for the system.
\label{thm:stab-mn}
\end{theorem}
 \begin{exmp} 
 	Consider the switching system $S$ of Example \ref{exmp:scalarswitching} with switching constrained by the graph of Figure \ref{fig:scalar12aut}, with nodes $a$ and $b$. 
 Here, the multinorm
  $$\mulno = \{ |x|_{a}, |x|_{b}  \} = 
  \left \{  4 |x| , 
    |x|
  \right \},
  $$  
   where $|x|$ is the absolute value of $x$, has a \emph{value} of $1/2$, obtained by applying (\ref{eq:mn}) to the system.
   Thus, the system is stable from Theorem \ref{thm:stab-mn}.
 \end{exmp}
 
\subsection{Extremal multinorms and boundedness}
Given a system  $S$, Proposition \ref{prop:CJSRmn} guarantees that the value of any multinorm is an upper bound on the CJSR of $S$. We now investigate the existence of multinorms with value \emph{equal} to the CJSR. The result proposed hereunder generalizes the characterization of \emph{extremal} norms for arbitrary switching systems (see e.g \cite{JuTJSR}, Section 2.1.2). 
\begin{theorem}
 A system $S(\graph(V,E), \mset)$ admits an \emph{extremal} multinorm, i.e. a multinorm 
 $\mulno^* = \{|\cdot|_v, v\in V\}$ with $\gamma^*(\mulno^*) = \jsr(S)$,  if and
 only if there exists a constant $K \geq 1 $ such that for any path $p$ in $\graph$,
 \begin{equation}
 |A_p x | \leq K \jsr(S)^{T_p}|x|,
 \label{eq:non-def}
 \end{equation}
 where $T_p \geq 0$ is the length of $p$.
 \label{thm:non-def}
 \end{theorem}
  \begin{pf}
  Assume first that $\jsr(S) > 0$. 
  By homogeneity of the CJSR with respect to the set $\mset$, we can further assume $\jsr(S) = 1$ without loss of generality (by scaling the matrices in $\mset$ by $1/\jsr(S) > 0$). \\
    We start with the only-if part. Let $\mulno^* = \{|{\cdot}|_v, v\in V\}$ be an \emph{extremal} multinorm.  Take any norm $|{\cdot}|$. By equivalence of norms in $\reels^n$, there are two scalars $\alpha, \beta > 0$ such that for any node $v \in V$, $\alpha |{\cdot}| \leq |{\cdot}|_v \leq \beta|{\cdot}| $. Since $\mulno$ is extremal, for any path $p$ between two nodes $v$ and $w$, we have $|A_p x|_w \leq  |x|_v$ from (\ref{eq:mn}), at which point we conclude that (\ref{eq:non-def}) holds for $K = \beta/\alpha$. \\*
For the if part, consider a system $S$ with $\jsr(S) = 1$ and define at each node $v \in V$ the following norm:
\begin{equation*}
|x|^*_v   = \sup_p \{|A_p x|: \, \, \text{$p$ is a path in $\graph$ starting at $v$}\},
\end{equation*}
where $|\cdot|$ is e.g. the Euclidean norm.
The functions $|{\cdot}|^*_v$, $v \in V$ are indeed norms  ($|x|^*_v \geq |x|$ since we also include paths of length 0).
Also, it is direct to check from the equation above that for any edge $(v, w, \sigma) \in E$, 
$|A_\sigma x|^*_w \leq |x|_v^*$.
Therefore, the multinorm $\mulno = \{ |{\cdot}|^*_v, \, v \in V \}$ is such that $\gamma^*(\mulno) \leq 1$. From Proposition  $\ref{prop:CJSRmn}$, since 
$\jsr(S) = 1$, we conclude that $\mulno$ is extremal.\\
In the case $\jsr(S) = 0$, we observe that (\ref{eq:non-def}) holds if and only if $A_p = 0$ for all path $p$ in $\graph$, which is equivalent to having the value of any multinorm equal to zero as well. Thus, we conclude the proof.
\end{pf} 
The result above provides a necessary  and sufficient condition for the \emph{boundedness} of constrained switching systems with $\jsr(S) = 1$,
 i.e. the existence of $K \geq 1$ such that for all $x \in \reels^n$ and path $p$ accepted by $\graph$, $|A_p x| \leq K |x|$. The boundedness of systems with $\jsr(S) = 1$ is known to be undecidable (see \cite{PhJuASCF} for sufficient conditions), and thus, so is the \emph{existence} of an extremal multinorm for a given switching system. 
In Subsection \ref{subset:suffext}, given a multinorm, we provide a sufficient condition for extremality, and apply it for computing the CJSR.

 \section{Approximation algorithms for stability analysis}
\label{sec:approx} 
Given a system $S(\graph, \mset)$ and a maximum \textit{relative} error $r > 0$, we wish to obtain an estimate $\tilde{\rho}$ of
the CJSR $\jsr(S)$, such that
$ \jsr(S) \leq \tilde{\rho} \leq (1+r) \jsr(S).$\\*
Through this section we will provide several approximation algorithms solving the above problem. 
All the methods share the same core mechanism: the approximation of a multinorm  for the system $S$, with value close to the CJSR, 
by a \textit{quadratic} multinorm $\mulno = \{|\cdot|_{Q, v}, v \in V\}$, where each norm is quadratic, i.e. $|x|_{Q,v} = x^\top Q_{v} x$
 for a positive definite matrix $Q_v \succ 0$. 
 This is  expressed as a \emph{quasi-convex} optimization program,
 solved by using a bisection procedure, iteratively checking the feasibility of the set of LMIs (\ref{eq:ptas_prog}). 
\begin{theorem}
Consider a system $S(\graph(V,E), \mset)$.\\ The value $\gamma_*(S)$ such that 
\begin{align}
&\gamma_*(S)  =  \inf_{\gamma, \{Q_v \in \reels^{n \times n}, v\in V\} }\gamma \nonumber\\
&\, s.t. \left \{ \begin{aligned}
  & \forall \, (v, w,  \sigma) \in E,  \\
 & \qquad -A_\sigma^\tran Q_w^{}A_\sigma^{} + \gamma^2 Q_v^{} \succeq 0, \\
 & \forall v\in V, \, Q_{v} \succ 0, \label{eq:ptas_prog} 
 \end{aligned} \right .
\end{align}
satisfies the following inequalities: 
\begin{equation}
 \hat \rho  (S) \leq \gamma_*(S) \leq \sqrt{n} \left (
 \hat \rho (S) \right ).
\label{eq:ptas1}
\end{equation}
Moreover, the LMI feasibility sub-problem 
(\ref{eq:ptas_prog}) is solved in a number of 
operations bounded by
\begin{equation}
\mathcal{O}\left ( n^{13/2}\cdot\card{V}^{7/2} \cdot \card{E}^{3/2} \right ).
\label{eq:complexity}
\end{equation}  
\label{thm:ptas1}
\end{theorem}
\vspace{-0.5cm}

\begin{pf}
First, we show that for any system $S(\graph, \mset)$ and any $\epsilon > 0$, 
there is a quadratic multinorm 
with a value $\gamma$ that satisfies
 $ \jsr(S) \leq \gamma \leq \sqrt{n} ( \jsr(S) + \epsilon).$\\
The result is obtained from 
\textit{John's Ellipsoid Theorem} \cite{JoEPWI}
 (see \cite{BlNeOTAO, PaJaAOTJ} for similar approaches 
 for arbitrary switching systems). 
 John's ellipsoid theorem states that for any norm $|{\cdot}|$ of $\reels^n$, 
 there exists a \emph{quadratic} norm
$|{\cdot}|_Q: x \rightarrow \left ( x^\tran Q x \right )^{1/2}$, with $Q \succ 0$, such that
$\forall x \in  \reels^n: \, |x|_Q \leq |x| \leq \sqrt{n} |x|_Q. $ 
Let us take $\epsilon > 0$ and a multinorm $\mulno_\epsilon = \{|{\cdot}|_v, v \in V\}$ with $\gamma^*(\mulno_\epsilon) \leq  \jsr(S) + \epsilon$.
Such a multinorm exists (Proposition \ref{prop:CJSRmn}).
By John's ellipsoid theorem,  there exist quadratic norms, forming a \textit{quadratic}
 multinorm $\mulno_{Q, \epsilon} = \{|{\cdot}|_{Q,v},  v \in V \}$, such that for any edge $(v, w, \sigma) \in E$, $\forall x \in  \reels^n:$
$$  |A_\sigma x|_{Q,w} \leq |A_\sigma x|_w \leq (\jsr(S) + \epsilon) |x|_v \leq \sqrt{n} (\jsr(S) + \epsilon) |x|_{Q,v}. $$
Since the above holds for any edge we can state that, for any $\epsilon > 0$,
 there is a quadratic multinorm $\mulno_{Q, \epsilon}$ with $ \gamma^*(\mulno_{Q, \epsilon}) \leq \sqrt{n} (\jsr(S) + \epsilon)$.
Taking $\epsilon \rightarrow 0$ we obtain (\ref{eq:ptas1}). \\*
The complexity computations are obtained from the reference book \cite{BeNeLOMC}, p.424.
The number of variables in the problem (\ref{eq:ptas1}) is $\left (n(n+1)/2 \right)\card{V}$. The LMIs constraints can be represented by a $ n \left ( \card{V}+\card{E} \right )$ bloc diagonal symmetric matrix with diagonal blocs of size $n \times n$.
\end{pf} \vspace{-0.5cm }
Notice that
1) the above result gives a sufficient condition for a given system to possess a quadratic \textit{Lyapunov} multinorm 
(if $\jsr(S) < 1/\sqrt{n}$, then (\ref{eq:ptas1}) guarantees $\gamma_*(S) < 1$)
 and that 2) it provides an algorithm to solve the relative approximation problem with maximum error $r \geq \sqrt{n} - 1$.
We now present ways to \emph{increase} the accuracy of the estimation, in Subsection \ref{subsec:Tlift} through Subsection \ref{subsec:pdep}, by performing an \textit{algebraic lifting} of the structures defining the system $S(\graph, \mset)$.
\subsection{ The $T$-product lift  }
\label{subsec:Tlift}
 This method allows to provide arbitrarily accurate estimates of the CJSR, with the cost of constructing a system on a graph with a large amount of edges.
  \begin{definition}[$T$-product lift]
  Given a system $S(\graph, \mset)$ and an integer $T \geq 1$, 
 the T-Product lift of $S$, denoted $S^T(\graph, \mset)$, is a constrained switching system on an automaton $\graph'(V',E')$ and a matrix set $\mset'$ defined as follows:
  \begin{enumerate}
  \item $\graph'$ has the same set of nodes as $\graph$ (i.e. $V' = V$).
  To each path $p$ in $\graph$ of length $T$, between two nodes $v$ and $w$ in $V$,
 is associated an edge $e = (v,w,\{\sigma(1) \ldots \sigma(T)\}) \in E'$. The label on this edge is a concatenation of those across the path $p$.
 \item The set of matrices $\mset'$ is the set of all products of size $T$ of matrices in $\mset$ that are accepted by $\graph$. For a label $\{\sigma(1) \ldots \sigma(T)\}$ of the lifted system, $A_{\{\sigma(1) \ldots \sigma(T)\}} = A_{\sigma(T)}\cdots A_{\sigma(1)} \in \mset'.$
  \end{enumerate}
  \label{def:tlift}
  \end{definition}  
If the system $S(\graph, \mset)$ describes the evolution of a state $x_t, t  =\{1, 2, \ldots\}$, then the system $S^T$ describes the evolution of the same state at times $kT$ for $k = \{1, 2, \ldots\}$. 
\begin{exmp}
The 2-product lift of the automaton of Figure \ref{fig:auto_ex1} is presented on 
Figure \ref{fig:auto_exTprod}.
\begin{figure}[!ht] █
\centering
\begin{subfigure}[c]{0.2\textwidth}
 \centering
\includegraphics[scale= 0.4]{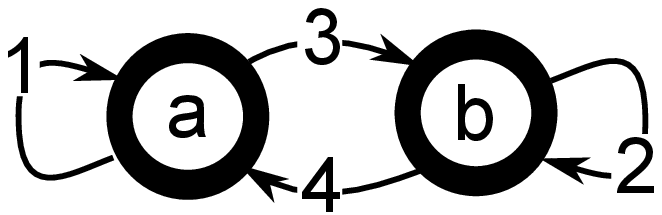}
\caption{}
\label{fig:auto_ex1}
\end{subfigure}
~
 \begin{subfigure}[c]{0.2\textwidth}
 \centering
\includegraphics[scale= 0.4]{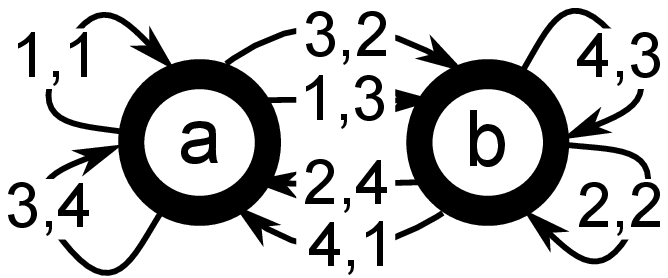}
\caption{}
\label{fig:auto_exTprod}
\end{subfigure}
\caption{The graph of Fig. \ref{fig:auto_exTprod} has one edge per path of length 2 in the graph of Fig. \ref{fig:auto_ex1}.}
\end{figure}
\end{exmp}

\begin{theorem}
The optimal value $\gamma_*(S^T)$ obtained by applying Theorem \ref{thm:ptas1} to the system $S^{T}( \graph, \mset)$ is such that
$ {\jsr(S) \leq \gamma_*^{1/T}(S^T) \leq  n  ^{1/(2T)} \jsr(S)}. $
\label{thm:prodbound}
\end{theorem}

\begin{pf}
We first show that for any system $S$ and any integer $T \geq 1$,
  $ \jsr(S^T) = \jsr(S)^T. $ Define for all $k \geq 1$, 
$$\jsr_k(S) = \max_{p} \{\|A_p\|^{1/k}: \, p \text{ of length $k$ accepted by $\graph$} \}.$$
Clearly, $\jsr_k(S^T) = \jsr_{kT}(S)^T,$ and by continuity of the exponentiation, 
we obtain $\jsr(S)^T = \lim_{k \rightarrow \infty}\jsr_{kT}(S)^T = \lim_{k \rightarrow \infty}\jsr_{k}(S^T) = \jsr(S^T) $.
Applying Theorem \ref{thm:ptas1} to $S^T$, the inequalities (\ref{eq:ptas1}) give
$ \jsr(S)^T = \jsr(S^T) \leq \gamma_*(S^T) \leq n^{1/2} \jsr(S^T) = n^{1/2} \jsr(S)^T,$
which concludes the proof.
\end{pf}

\begin{corollary}
For any system $S(\graph, \mset)$ and relative error bound $r > 0$, 
let
$$T = \lceil \log(n)/ (2 \log(1+r)) \rceil. $$
The CJSR estimation $\gamma_*(S^T)$ obtained by applying Theorem \ref{thm:ptas1}  to the system $S^{T}( \graph, \mset)$ satisfies
$ \jsr(S) \leq  \gamma_*(S^T)^{1/T} \leq (1+r) \jsr(S). $
 \label{cor:approxTot}
\end{corollary}
As a conclusion, given $r > 0$, there is a $T$-product lift of $S$ allowing the retrieve an estimate with relative error at most $1+r$. 
The amount of edges in the system $S^T$ increases exponentially with $T$, with one edge per path of length $T$. The amount of edges $\card{E}$ of  $\graph$ has however the least impact for computational  complexity (\ref{eq:complexity}).
\subsection{Approximation through sum-of-squares.}
\label{subsec:SOS}
In this subsection we generalize the technique presented in \cite{PaJaAOTJ} to constrained switching systems and devise a CJSR approximation scheme relying on an algebraic lifting of the set of matrices $\mset$. 
The procedure generating the lifted space is called the \emph{[$d$]-lift} (a full presentation can be found in \cite{PaJaAOTJ}, Section 3).  Given an integer $d$, the [$d$]-lift of a vector $x \in \reels^n$ is a vector $x^{[d]} \in \reels^{C(n+d-1, d)}$ of monomials of degree $d$, where ${C(k, \ell)}$   is the number of combinations of $\ell$ elements in a set of  $k$ elements. More important to us is that the \emph{[d]-lift} is well defined for linear maps, i.e. $\exists A^{[d]} : A^{[d]}x^{[d]} = (Ax)^{[d]}$,  this definition extends to sets of matrices $\mset^{[d]}$. Moreover,    $\left |x^{[d]} \right | = |x|^d$ for the Euclidean norm in the appropriate dimensions.\\  
In constrast with the method of Subsection \ref{subsec:Tlift},  we now conserve the graph of the constrained system, but approximate its CJSR by using potentially non-convex approximations of multinorms that are obtained from homogeneous sum-of-squares polynomials with degree $2d$,  $d \geq 1$ being an integer to be chosen. These polynomials have the form $x \mapsto (x^{[d]})^\top Qx^{[d]}$, where $Q \succ 0$ is a quadratic form in the lifted space.

\begin{theorem}
Given a system $S(\graph, \mset)$ and an integer $d \geq 1$, let $S^{[d]}(\graph, \mset)$ denote the constrained switching system on the same graph $\graph$ and on $[d]$-lift $\mset^{[d]}$ of $\mset$.\\*
The value $\gamma_*(S^{[d]})$ obtained by applying
 Theorem \ref{thm:ptas1} to the system $S^{[d]}$ satisfies
$ \jsr(S) \leq \gamma_*^{1/d}(S^{[d]}) \leq   C(n+d-1,d)  ^{1/(2d)} \jsr(S). $
 \label{thm:sosbound}
\end{theorem}
\begin{pf}
Since $|x^{[d]}| = |x|^{d}$ holds for the euclidean norm, we have $\jsr(S^{[d]}) = \jsr(S)^d$ from (\ref{eq:jsr}). Given the dimension of the set $\mset^{[d]}$, Theorem \ref{thm:ptas1} produces an estimate such that
$ \jsr(S^{[d]}) \leq \gamma_* \leq  C(n+d-1,d) ^{1/2} \jsr(S^{[d]}).$
\end{pf} 
The transformation here affects the dimension $n$ of the system. 
As it is shown in \cite{PaJaAOTJ} - Example 4, the CJSR approximation can be further refined by making use of explicit sum-of-square programming (rather than solving the program of Theorem \ref{thm:ptas1} in the lifted space), but with the same accuracy bounds.

 \subsection{Improving accuracy by adding memory to the graph.}
\label{subsec:pdep} 
Path-dependent Lyapunov functions have been introduced by Lee and Dullerud \cite{LeDuUSOD}
as tools for the stability analysis and control of discrete-time switching systems.
The concept, which follows a similar idea to that of Bliman and Ferrari-Trecate \cite{BlFeSAOD}, is to build a multiple Lyapunov function that associates a different quadratic form to each switching sequence of a length $M \geq 0$, $M$ being an integer parameter called the \textit{memory}  of the function.
The authors showed that, for any stable switching system, there is a finite $M$ such that the system admits a path-dependent Lyapunov function with memory $M$.\\*
Similar statements can be made about the T-Product lift and the $[d]$-lift defined in the previous subsections. Indeed, if $\jsr(S) < 1$ and since the approximations are asymptotically tight, there is a finite value $T$ (or $d$) for which the approximation algorithm will return an estimate lower then 1.  Given this parallel, it is natural to ask whether the methods of \cite{LeDuUSOD, BlFeSAOD} present similar approximation bounds to that of  Theorems  \ref{thm:prodbound} and \ref{thm:sosbound}.\\* 
We first define a lifting procedure, adapted from \cite{LeDuUSOD, BlFeSAOD}, allowing to obtain path-dependent Lyapunov function as \emph{quadratic multinorms of an augmented automaton}.  
\begin{definition}[M-Path-Dependent Lifting]
Given a system $S(\graph(V,E), \mset)$ and an integer $M \geq 0$, \emph{the $M$-Path-Dependent lift} of S, denoted $S_M(\graph, \mset)$, is a constrained switching system with the same set of matrices $\mset$ and with an automaton $\graph'(V',E')$ which is the same as $\graph$ for $M = 0$, and  constructed as follows for $M \geq 1$. Start with $V' = \{\}$, $E' = \{\}$, then:
\begin{enumerate}
\item For each path $p$ of length $M$ in $\graph$, add the node $v_p$ in $V'$.
\item For each path $p$ of length $M+1$ in $\graph$, let $p = (e_1, \ldots, e_{M+1})$, where $e_k$ is the $kth$ edge of the path. Define $p^- = (e_1, \ldots, e_M)$, $p^+ = (e_2, \ldots, e_{M+1})$, and $\sigma_p$ 
as the label of the edge $e_{M+1}$. Then, add the edge $(v_{p^-}, v_{p^+}, \sigma_p) \in E'$.
\end{enumerate}
\label{def:pdep}
  \end{definition}  
\begin{exmp}
Figure \ref{fig:auto_exPD} presents the 1-Path-Dependent lift $\graph'$ of the automaton $\graph$ of Figure \ref{fig:auto_ex1}. There are 4 nodes in $\graph'$, one for each edge of $\graph$. For the edge $(v, w, \sigma)$ of $\graph$, we use ``$vw$'' to refer to the corresponding node in $\graph'$.
For the edges, consider for example the path of length 2 $\left ( (b,b,2), (b,a,4)\right )$ in $\graph$. To this path corresponds the edge $(bb, ba, 4)$ in $\graph'$.  
 \begin{figure}[ht]
\centering
\includegraphics[scale= 0.4]{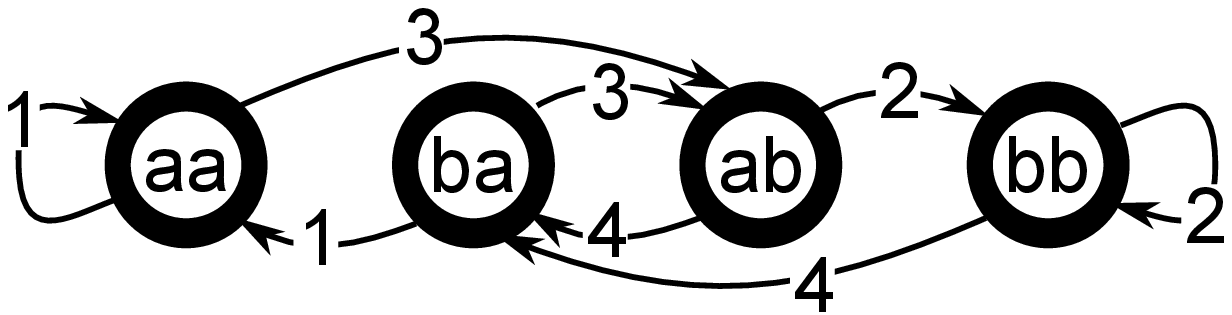}
\caption{}
\label{fig:auto_exPD}
\end{figure}
\end{exmp}

\begin{theorem}
Given a system $S(\graph, \mset)$ and an integer $T  \geq 1$, let $\gamma_{*}(S^T)$ and $\gamma_{*}(S_{T-1})$ be  the values obtained by applying Theorem \ref{thm:ptas1} to $S^T(\graph,\mset)$ and $S_{T-1}(\graph,\mset)$ respectively. The estimates satisfy $ \gamma_*(S_{T-1}) \leq \left ( \gamma_*(S^T) \right ) ^{1/T}. $
\label{thm:tliftandpdep}
\end{theorem}

\begin{pf}
For $T = 1$, $S = S_{T-1} = S^T$, so the claim holds. 
Assume now $T \geq 2$. 
We will show that given a  quadratic multinorm for $S^T$ with value $\gamma^T$, we can always construct a quadratic multinorm for $S_{T-1}$ with value at most $\gamma$. \\*
We refer to paths by their succession of edges, $p = (e_1, \ldots, e_T)$ is a path of length $T$ in $\graph$, and let $e_k = (v_{k-1}, v_{k}, \sigma_{k})$.
A quadratic multinorm with value $\gamma^T$ for $S^T$ associates a quadratic form $Q \succ 0$ per node in $\graph$ such that, for all paths $p$,
\begin{equation}
 A_p^\tran Q_{v_T}A_p - \gamma^{2T} Q_{v_0} \preceq 0.
 \label{eq:tliftlmi}
\end{equation}
To the same path $p$ corresponds an edge in the (T-1)-Path-Dependent lift, 
between a node $v_{ \left (e_1, \ldots, e_{T-1} \right)}$ and a node 
$v_{ \left (e_{2}, \ldots, e_{T} \right)}$.
For a quadratic multinorm of $S_{T-1}$ with value $\gamma$ to exist, there must be a quadratic form $Q_p \succ 0$ per path $p$ of length $T$, such that the following holds:
\begin{equation}
 \begin{aligned}
&A_{\sigma_T}^\tran Q^{}_{(e_2, \ldots, e_T)}A_{\sigma_T}^{} \\
&\qquad   {} - \gamma^{2} Q^{}_{ \left 
(e_{1}, \ldots, e_{T-1} \right )} \preceq 0.
\label{eq:pdep}
\end{aligned}
\end{equation}
Letting $R = Q^{-1}$ denote the inverses of the quadratic forms in the
corresponding lift, by a Schurr complement, the LMIs
(\ref{eq:tliftlmi}) are equivalent to
\begin{equation}
 A_pR_{v_0}A_p^\tran - \gamma^{2T} R_{v_T}, \preceq 0,
 \label{eq:tliftlmidu}
\end{equation}
and the LMIs (\ref{eq:pdep}) to
\begin{equation}
\begin{aligned}
&A_{\sigma_T} R^{}_{ \left (e_1, \ldots, e_{T-1} \right )}A_{\sigma_T}^\tran  \\
&\qquad   \qquad {} - \gamma^{2} R^{}_{ \left ( e_{2}, \ldots, e_T \right ) } \preceq 0.
\label{eq:pdepdu}
\end{aligned}
\end{equation}
Assume that we have a solution $\{R_v, v \in V\}$ to  (\ref{eq:tliftlmidu}).  Given any path $p = (e_1, \ldots, e_{T-1})$, with labels $\sigma_1, \ldots, \sigma_{T-1}$ that visits the nodes $v_0, \ldots, v_{T-1}$ in $\graph$,  define
\begin{equation}
\begin{aligned}
& R^{}_{ \left (e_1, \ldots, e_{T-1} \right )   }  =  R^{}_{v_{T-1}}  \nonumber \\
& \quad +  \gamma^{-2} A^{}_{ \sigma_{T-1} } R^{}_{v_{T-2}} A_{ \sigma_{T-1} }^\tran   \nonumber \\
& \quad {} +\gamma^{-4} A^{}_{\left ( \sigma_{T-2}, \sigma_{T-1} \right )} 
  R^{}_{v_{T-3}} 
   A^{\tran}_{\left ( \sigma_{T-2}, \sigma_{T-1} \right )} \nonumber \\
& \quad {} + \cdots\\
& \quad {} +\gamma^{-2(T-1)}  A^{}_{\left ( \sigma_1, \ldots, \sigma_{T-1} \right )}  R^{}_{v_0}  A^{\tran}_{\left ( \sigma_1, \ldots, \sigma_{T-1} \right )},
\label{eq:lmifin}
\end{aligned}
\end{equation}
where $A_{(\sigma_k, \ldots, \sigma_{T-1})} = A_{\sigma_{T-1}} \cdots A_{\sigma_{k}}$.
Injecting these quadratic forms in (\ref{eq:pdepdu}), we obtain (\ref{eq:tliftlmidu}), and thus we have a quadratic multinorm for $S_{T-1}$, with value at most $\gamma$. 
\end{pf}
\begin{corollary}
Given a system $S(\graph, \mset)$ and an integer $M \geq 1$, the value $\gamma_*$ obtained by applying Theorem \ref{thm:ptas1} to the $M$-Path-Dependent lift of S satisfies
$ \jsr(S) \leq \gamma_*  \leq n^{1/(2(M+1))} \jsr(S).$
\label{cor:pdepandtlift}
\end{corollary}
The amount of nodes $\card{V}$ and edges $\card{E}$ in the M-Path-Dependent lift both grow exponentially with $M$, with one node per path of length $M$, and one edge per path of length $M+1$ in $\graph$.  
\subsection{A sufficient condition for extremality of quadratic multinorms.}
\label{subset:suffext}
We now present an easy to check sufficient condition under which the CJSR estimate obtained by applying Theorem \ref{thm:ptas1} is exact.
The condition can be applied for all the lifts developed above, since they all rely on Theorem \ref{thm:ptas1} to retrieve a quadratic multinorm of minimal value. 
We assume that the estimate of Theorem \ref{thm:ptas1} is attained by the value of a quadratic multinorm, and give a sufficient condition for its extremality.
We start with the following observation.
\begin{lemma}
For any system $S(\graph, \mset)$ and any cycle\footnote{A cycle is a closed path, i.e. whose source and destination nodes are the same. } $c$ of length $T$ in $\graph$, $ \jsr(S) \geq \rho(A_c)^{1/T}$,
where $\rho(A_c)$ is the spectral radius of the product $A_c$.
\label{lemma:cycle}
\end{lemma}
\begin{pf}
For any induced matrix norm $\|\cdot\|$, the following holds
$$
\begin{aligned}
 \rho(A_c)^{1/T} &= \lim_{k \rightarrow \infty} \|A_c^k\| ^{1/kT} \\
& \leq \lim_{k \rightarrow \infty} \max\{\|A_p\|^{1/kT}: p \text{ is a path of length $kT$}\},\\
& \leq \jsr(S).
\end{aligned} 
$$
\end{pf}

We now define a \emph{simple} cycle in  the graph $\graph$ as a cycle such that for every node $v$ visited by the cycle,
there exists no partition of the cycle into two cycles on $v$.

\begin{theorem}[Sufficient extremality condition]
Let $\{Q_v, \, v \in V\}$ be the optimal quadratic forms obtained by applying Theorem \ref{thm:ptas1} to the system $S(\graph(V,E), \mset)$, corresponding to a multinorm $\mulno$ with value $\gamma_*$. If the set of edges 
$$ E' = \left \{(v, w, \sigma) \in E: \lambda_{\min}\left ( -A_\sigma^\tran Q_w^{}A_\sigma^{} + \gamma^{2}_* Q_v^{} \right ) = 0 \right \},$$
 where $\lambda_{\min}(X)$ denotes the smallest eigenvalue of the positive definite matrix $X$,
 forms a \emph{simple} cycle $c$ in $\graph$, 
then $\mulno$ is \emph{extremal}, and $\jsr(S) = \rho(A_c)^{1/T},$
where $T$ is the length of $c$.
\label{prop:sufcycles}
\end{theorem}

\begin{pf}
The edges in the set $E'$ are  those for which the LMI constraints in (\ref{eq:ptas_prog}) are tight for the multinorm $\mulno$. 
It is a known fact of convex optimization that removing constraints that \emph{are not} tight at the optimal solution of a given program
\emph{does not} affect its optimal solutions. 
Therefore, it must be the case that the optimal values $\gamma_*$ and $\gamma_*'$ obtained  by applying Theorem \ref{thm:ptas1} respectively to $S(\graph(V,E), \mset)$ and $S(\graph(V,E'), \mset)$ are \emph{equal}, and we may focus on the second system. \\*
Since  $\graph(V,E')$ is a simple cycle, then for any $M \geq 0$, the M-Path-Dependent lift (Definition \ref{def:pdep}) leaves the graph invariant. Applying Corollary \ref{cor:pdepandtlift} with $M \rightarrow \infty$, we conclude that the multinorm obtained from applying Theorem \ref{thm:ptas1} to $S(\graph(V,E'), \mset)$ needs to be extremal for this cyclic graph, and $\gamma_*' = \rho(A_c)^{1/T}$.\\
Having computed $\gamma_*'$, and since $\gamma_* = \gamma'_*$, we now consider the original system. Applying Lemma \ref{lemma:cycle} and Proposition \ref{prop:CJSRmn}, we get
$ \rho(A_c)^{1/T}  \leq \hat \rho(\graph, \mset) \leq \gamma_* = \rho(A_c)^{1/T},$
which concludes the proof.
\end{pf} 

\subsection{The Kroenecker lift}
\label{subset:arblift}
We end this section by presenting another approach to the stability analysis of constrained switching systems.  Kozyakin \cite{KoTBWF} and Wang \cite{WaRoSOLA} independently introduced a lifting procedure that allows to obtain, from any system $S$, a set of matrices $\mset_S$ such that $\jsr(S) = \jsr(\mset_S)$.   
\begin{definition}[Kroenecker lift]
 Given a system $S(\graph(V,E), \mset)$, with $\card{V}$ nodes $\{v_i \in V,\, 1 \leq i \leq \card{V}\}$ and $\mset \subset \reels^{n\times n}$,
the \emph{Kroenecker lift} of $S$ is a matrix set defined as 
$\mset_{\graph} = \{A_{(v_i,v_j,\sigma)}, \, (v_i,v_j,\sigma) \in E\},$
with for the edge $(v_i, v_j, \sigma)$, 
$$A_{(v_i, v_j, \sigma)} = \left ( \mathbf{e}(j)_{}^{} \mathbf{e}(i)^{\tran}\right) \kron A_{\sigma},$$
where $\mathbf{e}(k) \in \reels^{\card{V}}$ is the $k$th vector of the canonical basis of $\reels^{\card{V}}$, and $\kron$ is the Kroenecker product.
\label{def:lift}
\end{definition}
The methods in Subsection \ref{subsec:Tlift}, \ref{subsec:SOS} and \ref{subsec:pdep} are based on the ones of \cite{BlNeOTAO, LeDuUSOD, PaJaAOTJ} and make use of the concept of multinorm for providing estimates of $\jsr(S)$. It is natural to compare these methods with a direct application of the ones of \cite{BlNeOTAO, LeDuUSOD, PaJaAOTJ} for approximating $\jsr(\mset_\graph)$.
We will focus on \cite{BlNeOTAO}, which approximates the JSR of a set of matrices using Theorem \ref{thm:ptas1} with a single quadratic form. The conclusions naturally carry on to the other methods of \cite{PaJaAOTJ,LeDuUSOD}.
\begin{proposition}
Consider a system $S(\graph(V,E), \mset)$ and its associated set $\mset_\graph$. There exists $Q_\graph  \succ 0$ and  $\gamma > 0$ such that
$\forall A_{(v,w,\sigma)} \in \mset_S$,
$$ A_{(v,w,\sigma)}^\tran Q_\graph A_{(v,w,\sigma)} - \gamma^2 Q_\graph \preceq 0,$$
if and only if $S$ has a quadratic multinorm with value at most $\gamma$.
\label{prop:liftisequal}
\end{proposition}\vspace{-0.4cm}
 \begin{pf} 
 Assume without loss of generality that $\gamma = 1$. To ease the reading, we assume $n = 1$, so that for $e \in E$, $A_e \in \reels^{\card{V} \times \card{V}}$. For $e = (v,w,\sigma) \in E$, we let $B_e = A_e^\top Q_\graph^{} A_e^{}  - Q_\graph$. The element at row $k$ and column $\ell$  of $A_e$ is $A_e^{[k,\ell]}$, and $B_e^{[k, \ell]}$, $Q_\graph^{[k, \ell]}$ are defined similarly.
 Define the set of nodes as $\{v_i \in V, \, 1 \leq i \leq \card{V}\}$. For any $1 \leq k, \ell \leq \card{V}$, we have
  \begin{equation*}
 \begin{aligned}
  B_{(v_i, v_j, \sigma)}^{ [k, \ell]} 
&  = \sum_{i',j'}A_{(v_i, v_j, \sigma)}^ {\tran, [k, i']}
  Q_{{\graph}}^{[ i',j' ]} A_{(v_i, v_j, \sigma)}^ { [j', \ell]}
 - Q_{{\graph}}^{[k,\ell]}, \\
 & = \delta_{k,i} \delta_{ \ell, i} \left (A_{(v_i, v_j, \sigma)}^ {\tran, [k, j]}
  Q_{{\graph}}^{[ j,j ]} A_{(v_i, v_j, \sigma)}^ { [j, \ell]} \right)
 - Q_{{\graph}}^{[k,\ell]},
 \end{aligned}
\end{equation*} 
where $\delta_{i,j} = 1$ if and only if $i = j$ (Kroenecker delta). \\*
For the only if part, $Q_\graph \succ 0$ and $B_e \preceq 0$. 
Since $Q_\graph \succ 0$ and $B_e \preceq 0$, then for all   $1 \leq k \leq \card{V}$, $Q_\graph^{[k,k]} > 0$, and
for any $e = (v_i, v_j, \sigma) \in E$, 
$B^{[i,i]}_e = A^{\tran}_\sigma
Q_\graph^{[j,j]} A_\sigma - Q_\graph^{[i,i]} \leq 0.$ Thus, we extract a quadratic multinorm for $S$ using the quadratic forms $Q_\graph^{[i,i]}$, $1 \leq i \leq \card{V}$\\*
For the if part of the proof, we take quadratic multinorm with value lower than 1. Let 
$Q_i > 0, 1 \leq i \leq \card{V}$, be the quadratic forms associated at each node. We reconstruct $Q_\graph$ by setting $Q_\graph^{[k,\ell]} = Q_{k}$ if $k = \ell$, and   $Q_\graph^{[k,\ell]} = 0$ else.
 \end{pf}  
 Applying Theorem 3.1 to $\mset_\graph$ therefore produces the same estimate as if it was applied on $S$. However, when using $\mset_\graph \subset \reels^{n \card{V} \times n \card{V}}$, the  complexity (\ref{eq:complexity}) reads $\bigO \left ( (n\card{V})^{(13/2)}  \card{E}^{(3/2)} \right )$, a significant increase compared to the previous $\bigO \left ( n^{(13/2)} \card{V}^{(7/2)} \card{E}^{(3/2)} \right )$. Moreover, the accuracy bounds (\ref{eq:ptas1}) become worse, due to the larger dimension of $\mset_\graph$.

\section{Numerical example}
\label{sec:example}
\begin{figure}[!ht]
\centering
\includegraphics[scale = 0.4]{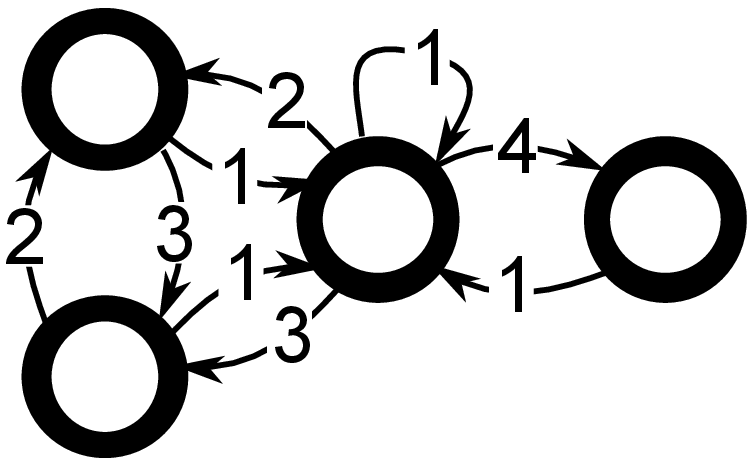}
\caption{}
\label{fig:dropout}
\end{figure}

\begin{figure}[!ht]
\centering
\includegraphics[scale = 0.45]{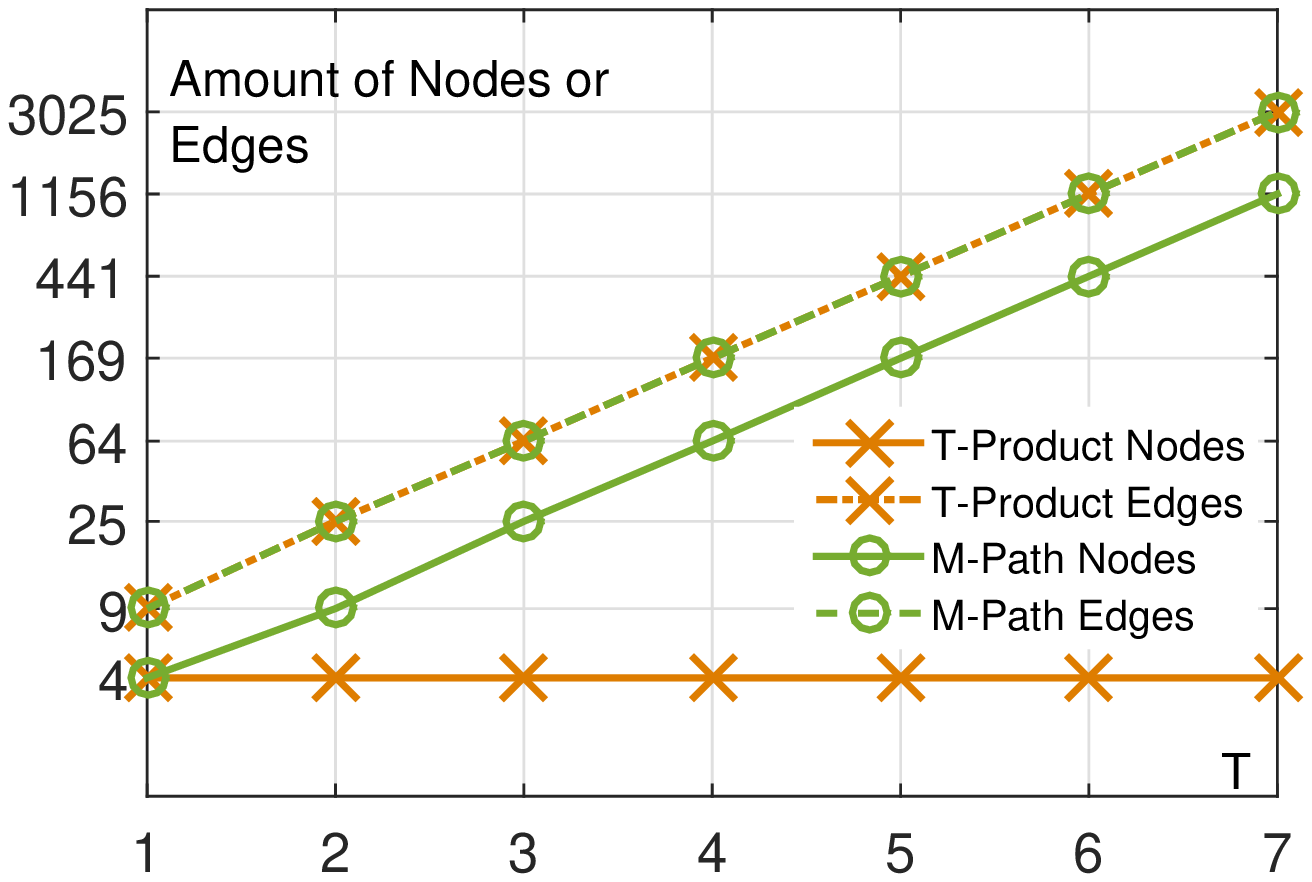}
\caption{}
\label{fig:compSize}
\end{figure}

\begin{figure}[!ht]
\centering
\includegraphics[scale = 0.45]{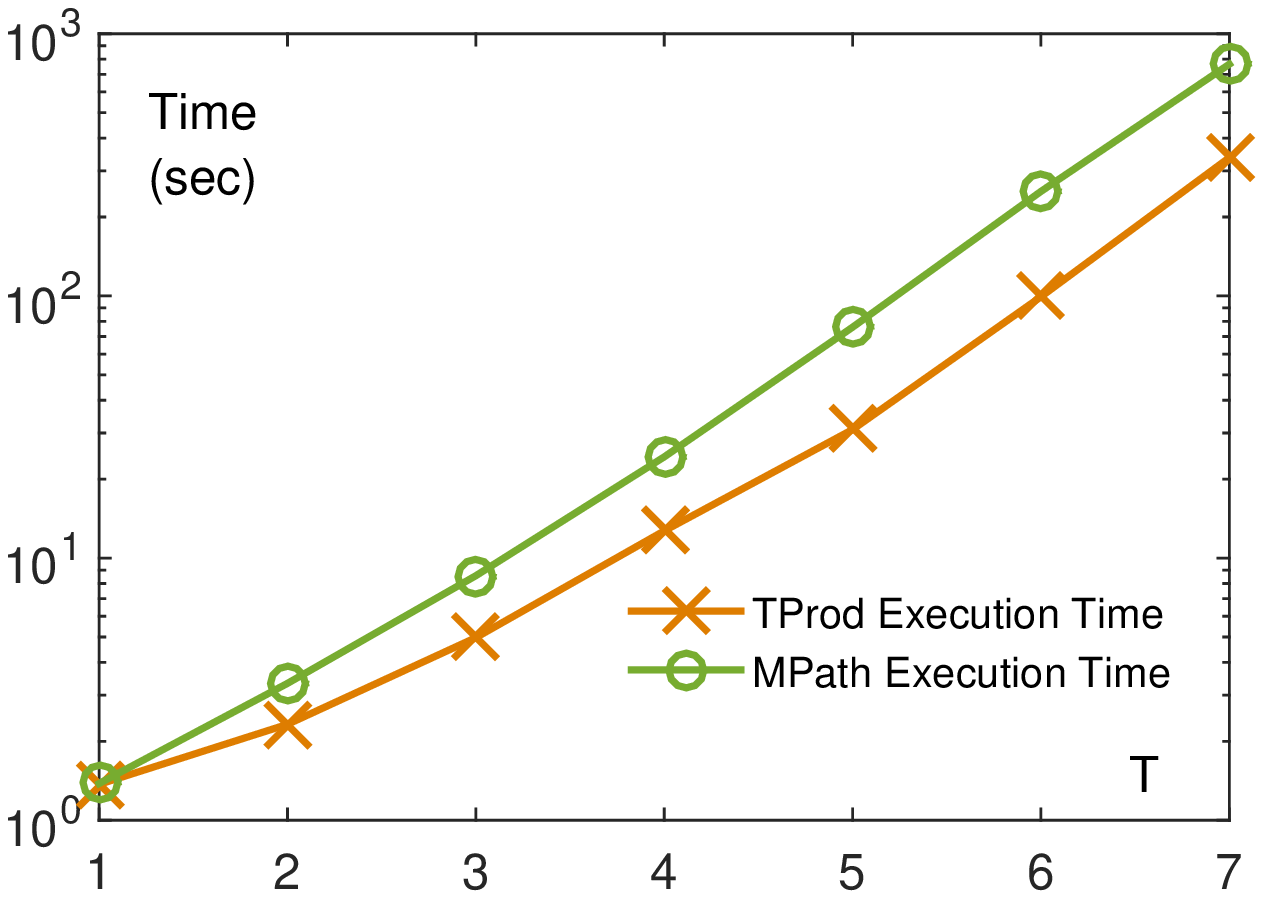}
\caption{}
\label{fig:compTime}
\end{figure}

\begin{figure}[!ht]
\centering
\includegraphics[scale = 0.45]{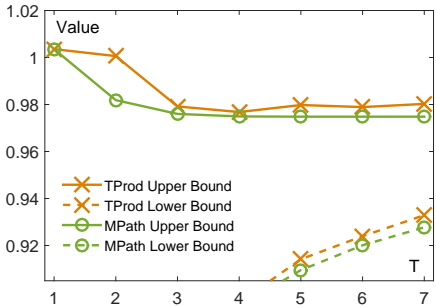}
\caption{Using the M-Path-Dependent lifts, we produce better estimates. The estimates produced by the T-Product lifts need not decrease monotonically as $T$ increase. Lower bounds computed from the relative accuracies of the CJSR estimates.}
\label{fig:compEstim}
\end{figure}

We consider the dynamics of a plant that may experience controller failures: $x_{t+1} = \left ( A + BK_{\sigma(t)} \right ) x_t,$
with 
$$ A = \begin{pmatrix} 0.94 & 0.56 \\ 0.14 & 0.46 \end{pmatrix}, \, B = \begin{pmatrix} 0 \\ 1 \end{pmatrix},$$
and $K_{\sigma(t)} = (k_{1, \sigma(t)}, k_{2, \sigma(t)} )$. The control gains switch to represent 4 failures modes. 
When everything works as expected, $\sigma(t) = 1$, and $$K_1 =  \begin{pmatrix} k_1 & k_2 \end{pmatrix} = 
\begin{pmatrix}
-0.49 & 0.27
\end{pmatrix}.  $$
The second and third modes correspond respectively to a failure of the first and second part of the controller, with $K_2 = (
0,  k_2
)$ and $K_3 = (
k_1,  0
).$  The last mode is the total failure case with $K_4 = (
0,  0
).$ 
As a constraint, we consider that a same part of the controller never fails twice in a row. We  obtain a constrained switching system on the set 
$\mset = \{A+BK_\sigma\}, \, \sigma  \in \{ 1, \ldots, 4\},$
with the automaton $\graph(V,E)$ depicted in Figure \ref{fig:dropout}.\\
We first  provide increasingly accurate estimates of the CJSR of the system, and then  provide the exact value of the CJSR  using Theorem \ref{prop:sufcycles}.
For the estimations, we compare  the T-Product lift of Subsection \ref{subsec:Tlift}, and the M-Path-Dependent lift, introduced in \cite{LeDuUSOD} and further studied in Subsection  \ref{subsec:pdep}. We apply these methods for $T = 1,\ldots,7$, and $M = T - 1 = 0, \ldots, 6$. From Corollary \ref{cor:approxTot} and Theorem \ref{cor:pdepandtlift}, these choices produce estimates within  $[\jsr(S), (1+r)\jsr(S)]$, with $1+r$ ranging from $\sqrt{2} \simeq 1.41$ for $T = 1$ to $2^{1/14} \simeq 1.05$ for $T = 7$. \\
For each value of $T$ and $M$ we first compute the lifted systems $S^T$ and $S_M$. Figure \ref{fig:compSize} provides a comparison of the amount of nodes and edges of these systems. Then, we solve the optimization program of Theorem \ref{thm:ptas1} on the lifts, applying Theorems \ref{thm:prodbound} and \ref{thm:tliftandpdep} to obtain bounds on the CJSR.\footnote{Matlab codes reproducing the results available at\\``\textbf{\url{http://sites.uclouvain.be/scsse/postrev.zip}}''.} The execution times for producing the estimates are compared in Figure \ref{fig:compTime}, and the estimation values are compared in Figure \ref{fig:compEstim}. The best estimation is achieved with the $M$-Path-Dependent lift, for $M = 6$, guaranteeing that
$ \jsr(S) \leq 0.9748\ldots.$ However, the T-Product lift provides estimates much faster.\\
We now  use Theorem \ref{prop:sufcycles} for computing the exact value of the CJSR. The results of Figure \ref{fig:compEstim} indicate that the sufficient condition can  only be met by the M-Path-Dependent lift, for $M = T-1 = 4,5$ or $6$.
In this case, for $M = 5$, the conditions of Theorem \ref{prop:sufcycles} are met. The obtained simple cycle corresponds to the switching sequence repeating the labels
$\{2,3,1,1,1,1,2,1\} $ and reaches the CJSR. 
Indeed, we have (Lemma \ref{lemma:cycle}) 
$$ \begin{aligned}
 \jsr(S) & \geq \rho(A_1A_2A_1^4A_3A_2)^{1/8} \simeq 0.9478\ldots,
\end{aligned}$$
which matches  our best CJSR estimate.\\
\section{Conclusion}
We  established a new framework for the stability analysis of discrete-time linear switching systems with constrained switching sequences.
It relies around the newly introduced concept of \emph{multinorms} and their link to the constrained joint spectral radius \cite{DaAGTS} (CJSR).
 The stability of constrained switching systems is equivalent to the existence of a set of exactly \emph{one norm per node} of the automaton, 
 with contractivity relations that are given by the edges of the automaton.\\*
By approximating these norms individually with quadratic norms,
 we design  the first arbitrarily accurate CJSR approximation schemes. 
 The framework also encapsulate well-known methods such as \emph{path-dependent} Lyapunov functions \cite{LeDuUSOD, BlFeSAOD} as methods approximating  extremal multinorms. \\*
In the future, we will apply this framework to control-oriented problems.
 Path-dependent Lyapunov functions have been used  in feedback controller design \cite{LeDuUSOD, EsLeCOLS},  and our goal is now to  give guarantees on the  performance of these controllers. Their usage for the study of the \emph{contractiveness} of switching  systems \cite{LeDuODAF, EsLeCOLS} suggests an estimation framework in this setting.

\bibliographystyle{IEEEtranS}
\bibliography{biblio}
\section{Annex: Proof of Theorem \ref{thm:stability-cjsr}}
\label{annex:proofCJSR}
From Definition 2 in \cite{DaAGTS}, we know that the limit (\ref{eq:jsr}) converges.
 Since the limit (\ref{eq:jsr}) converges, 
$ \forall \epsilon > 0, \, \exists T_\epsilon \geq 0: \, \forall t \geq T_\epsilon$, 
\begin{equation*}
\jsr(S)^t \leq \max_{\sigma({\cdot}) \, \text{accepted by $\graph$}} \|A_{\sigma(t-1)} \cdots A_{\sigma(0)}\| \leq (\jsr(S) + \epsilon)^t.
\end{equation*}
To show that $\jsr(S) < 1$ implies (exponential) stability, it suffice to take $\epsilon < 1 - \jsr(S)$. Indeed, we then  get that for all accepted sequences, $\lim_{t \rightarrow \infty}\|A_{\sigma(t-1)}\cdots A_{\sigma(0)}\| = 0.$ Exponential stability is then acquired since there can only be a finite amount of products of length $t \leq T_\epsilon$ with $\|A_{\sigma(t-1)} \ldots A_{\sigma(t)}\| > (\jsr(S) + \epsilon)^t$.
\\
Consider now the case $\jsr(S) \geq 1$.
For all $t \geq 1$, we define $x_{*}^t$
as 
$$x_{*}^t = \text{arg}\max_{|x| = 1} \max_{\sigma({\cdot}) \, \text{accepted by $\graph$}}
|A_{\sigma(t-1)}\cdots A_{\sigma(0)}x|.$$
We extract from the sequence $\{x_{*}^t\}_{t \geq 1}$ a subsequence converging to a point $x_*$, $|x_*| = 1$. From this point, there is a switching sequence satisfying  
$\lim_{t \rightarrow \infty}|A_{\sigma(t-1)}\cdots A_{\sigma(0)} x_*| \geq 1/2.$ Thus, the system is not asymptotically stable when $\jsr(S) \geq 1$.\\
This concludes the proof of Theorem $\ref{thm:stability-cjsr}$.

\end{document}

%% file: pkg_hdr.tex
\usepackage{amsmath}									
\usepackage{amssymb}
\usepackage{cite}
\usepackage[usenames,dvipsnames]{color}
\usepackage{pgfplots}
	\pgfplotsset{compat=newest} 
	\pgfplotsset{plot coordinates/math parser=false}
\usepackage{pstool} 
\usepackage{psfrag}
\usepackage{graphicx}
\usepackage{epstopdf} 
	\DeclareGraphicsExtensions{.pdf,.png,.jpg,.eps,.ps}

\usepackage{comment}
\usepackage{caption}
\usepackage{subcaption}
\usepackage[hidelinks]{hyperref}